\documentclass[11pt,a4paper,oneside]{article}

\usepackage[latin1]{inputenc}
\usepackage{subfigure}
\usepackage{multicol}

\usepackage{graphicx}
\usepackage{amssymb}

\usepackage{amsmath, xypic}

\newtheorem{thm}{Theorem}

\newtheorem{exm}{Example}

\newtheorem{rem}{Remark}

\newcommand{\al}{\alpha}

\newcommand{\la}{\lambda}

\newcommand\GenBin[2]{\left \langle \begin{array}{cc}#1\\#2\end{array} \right \rangle}
\newcommand\VarGenBin[2]{\left [ \begin{array}{cc}#1\\#2\end{array} \right ]}
\newlength{\cellsz}
\newcounter{cellsize}
\newcommand{\setcellsize}[1]{%
  \setcounter{cellsize}{#1}%
  \setlength{\cellsz}{\value{cellsize}\unitlength}}%
\setcellsize{16}%
\newcommand\cellify[1]{\def\thearg{#1}\def\nothing{}%
\hbox to 0pt{{\begin{picture}(\value{cellsize},\value{cellsize})
  \put(0,0){\line(1,0){\value{cellsize}}}
  \put(0,0){\line(0,1){\value{cellsize}}}
  \put(\value{cellsize},0){\line(0,1){\value{cellsize}}}
  \put(0,\value{cellsize}){\line(1,0){\value{cellsize}}} \end{picture}
\hss}}
\vbox to \cellsz{ \vss \hbox to \cellsz{\hss$#1$\hss} \vss}}
\newcommand\tableau[1]{\vcenter{\vbox{\let\\\cr
\baselineskip -16000pt \lineskiplimit 16000pt \lineskip 0pt
\ialign{&\cellify{##}\cr#1\crcr}}}}
\newcommand\tabl[1]{\vtop{\let\\\cr
\baselineskip -16000pt \lineskiplimit 16000pt \lineskip 0pt
\ialign{&\cellify{##}\cr#1\crcr}}}
\newlength{\varcellsz}
\newcounter{varcellsize}
\newcommand{\setvarcellsize}[1]{%
  \setcounter{varcellsize}{#1}%
  \setlength{\varcellsz}{\value{varcellsize}\unitlength}}%
\setvarcellsize{10}%
\newcommand\varcellify[1]{\def\varthearg{#1}\def\varnothing{}%
\hbox to 0pt{{\begin{picture}(\value{varcellsize},\value{varcellsize})
  \put(0,0){\line(1,0){\value{varcellsize}}}
  \put(0,0){\line(0,1){\value{varcellsize}}}
  \put(\value{varcellsize},0){\line(0,1){\value{varcellsize}}}
  \put(0,\value{varcellsize}){\line(1,0){\value{varcellsize}}} \end{picture}
\hss}}
\vbox to \varcellsz{ \vss \hbox to \varcellsz{\hss$#1$\hss} \vss}}
\newcommand\vartableau[1]{\vcenter{\vbox{\let\\\cr
\baselineskip -16000pt \lineskiplimit 16000pt \lineskip 0pt
\ialign{&\varcellify{##}\cr#1\crcr}}}}
\newcommand\vartabl[1]{\vtop{\let\\\cr
\baselineskip -16000pt \lineskiplimit 16000pt \lineskip 0pt
\ialign{&\varcellify{##}\cr#1\crcr}}}

\author{Ekaterina A. Vassilieva}
\title{Explicit monomial expansions of the generating series for connection coefficients}

\begin{document}
\maketitle
\begin{abstract}
This paper is devoted to the explicit computation of generating series for the connection coefficients of two commutative subalgebras of the group algebra of the symmetric group, the class algebra and the double coset algebra.
As shown by Hanlon, Stanley and Stembridge (1992), these series gives the spectral distribution of some random matrices that are of interest to statisticians. Morales and Vassilieva (2009, 2011) found explicit formulas for these generating series in terms of monomial symmetric functions by introducing a bijection between partitioned hypermaps on (locally) orientable surfaces and some decorated forests and trees. Thanks to purely algebraic means, we recover the formula for the class algebra and provide a new simpler formula for the double coset algebra. As a salient ingredient, we derive a new explicit expression for zonal polynomials indexed by partitions of type $[a,b,1^{n-a-b}]$. 
\end{abstract}

\section{Introduction}
\label{sec : intro}
For integer $n$ we note $S_n$ the symmetric group on $n$ elements and $\lambda = (\lambda_1, \lambda_{2},...,\lambda_p) \vdash n$ an integer partition of $\ell(\lambda) = p$ parts sorted in decreasing order. We define as well $Aut_\lambda = \prod_i m_i(\lambda)!$ and $z_\lambda =\prod_i i^{m_i(\lambda)}m_i(\lambda)!$  where $m_i(\lambda)$ is the number of parts in $\lambda$ equal to $i$. Let $m_\lambda(x)$ be the monomial symmetric function indexed by $\lambda$ on indeterminate $x$, $p_\lambda(x)$ and $s_\lambda(x)$  the power sum and schur symmetric function respectively.
Let $C_\lambda$ be the conjugacy class of $S_n$ containing the permutations of cycle type $\lambda$. The cardinality of the conjugacy classes is given by $|C_\la| = n!/z_\la$. Additionally, $B_n$ is the hyperoctahedral group (i.e the centralizer of $f_\star = (12)(34)\ldots(2n-1\,2n)$ in $S_{2n}$). We note $K_\lambda$ the double coset of $B_n$ in $S_{2n}$ consisting in all the permutations $\omega$ of $S_{2n}$ such that $f_\star \circ\omega\circ f_\star\circ\omega^{-1}$ has cycle type $(\lambda_1,\lambda_1, \lambda_{2},\lambda_{2},...,\lambda_p,\lambda_p)$. We have $|B_n| = 2^nn!$ and $|K_\la| = |B_n|^2/(2^{\ell(\la)}z_\la)$.
By abuse of notation, let $C_\lambda$ (resp. $K_\lambda$) also represent the formal sum of its elements in the group algebra $\mathbb{C} S_{n}$ (resp. $\mathbb{C} S_{2n}$). Then, $\{C_\lambda, \lambda \vdash n\}$ (resp. $\{K_\lambda, \lambda \vdash n\}$) forms a basis of the class algebra (resp.  double coset algebra).
For $\lambda$, $\mu$, $\nu \vdash n$, we focus on the connection coefficients $c^\nu_{\lambda\mu}$ and $b^\nu_{\lambda\mu}$ than can be defined formally by:
\begin{equation}
c^\nu_{\lambda\mu} = [C_\nu]C_\lambda C_\mu, \;\;\;\;\; b^\nu_{\lambda\mu} = [K_\nu]K_\lambda K_\mu
\end{equation}
From a combinatorial point of view $c^\nu_{\lambda\mu}$ is the number of ways to write a given permutation $\gamma_\nu$ of $C_\nu$ as the ordered product of two permutations $\alpha\circ\beta$ where $\alpha$ is in $C_\lambda$ and $\beta$ is in $C_\mu$. Similarly, $b^\nu_{\lambda\mu}$ counts the number of ordered factorizations of a given element in $K_\nu$ into two permutations of $K_\lambda$ and $K_\mu$. 
Although very intractable in the general case, Morales and Vassilieva (\cite{MV09}) found an explicit expression of the generating series for the $c^\nu_{\lambda\mu}$ for the special case $\nu = (n)$ in terms of monomial symmetric functions.

\begin{thm}[Morales and Vassilieva, 2009]
\label{thm : mv09}
\begin{equation}
\frac{1} {n}\sum_{\la, \mu \vdash n} {c_{\la\mu}^n} p_{\la}(x)p_{\mu}(y) =  \sum_{\la, \mu \vdash n} \frac{(n-\ell(\la))!(n-\ell(\mu))!}{(n+1-\ell(\la)-\ell(\mu))!}m_{\la}(x)m_{\mu}(y)
\end{equation}

\end{thm}

In \cite{FV10}, Feray and Vassilieva found an interesting expression when we set $\mu = (n)$ as well :

\begin{thm}[Feray and Vassilieva, 2010]
\label{thm : fv10}
\begin{equation}
\frac{1} {n!}\sum_{\la \vdash n} {c_{\la,n}^n} p_{\la}(x) =  \sum_{\la\vdash n} \frac{m_{\la}(x)}{n+1-\ell(\la)}
\end{equation}
\end{thm}

Both of these works rely on purely bijective arguments involving the theory of hypermaps on orientable surfaces. As shown in section \ref{sec : class}, these formulas can be recovered in a very simple way using some known relations between the connection coefficients, the characters of the symmetric group and the classical bases of the symmetric functions.\\
\indent The evaluation of the generating series for the $b^\nu_{\lambda\mu}$ is much more complicated. Morales and Vassilieva (in \cite{MV11}) found the first explicit formula in terms of monomial symmetric functions for the case $\nu = (n)$ thanks to a bijection between partitioned hypermaps in locally orientable surfaces and some decorated forests. 
Using an algebraic method analog to the class algebra case and a new formula for zonal polynomials on near hooks, section \ref{sec : doublecoset}, \ref{sec : zon poly} and \ref{sec : results} prove a new much simpler expression. To state it, we need to introduce a few more notations. If $x$ and $y$ are non negative integers, we define:
\begin{equation}
\GenBin{x}{y} = \frac{\binom{x}{y}^2}{\binom{2x}{2y}}, \;\;\;\;\; \VarGenBin{x}{y} = \frac{\binom{x}{y}^2}{\binom{2x+1}{2y}}
\end{equation}
\noindent if $x \geq y$ and $0$ otherwise. We also define the rational functions:
\begin{align}
&R(x,y,z,t,w) = \frac{(2x+w)(2y+w)(2z+w-1)(2t+w-1)}{(2x+w-1)(2y+w+1)(2z+w-2)(2t+w)}\\
&r_n(x,y) = 2n\frac{(n+x-y+1)(n+y-x)(n-x-y)!(2x-1)!!(2y-2)!!}{(-1)^{n+1-x-y}(n+x-y)(n+y-x-1)(2(x-y)+1)}
\end{align}

\noindent where the later definition of $r_n(x,y)$ is valid for $y\geq1$. Additionally, we define $r_n(n,0) = (2n-1)!!$. We have:

\begin{thm}[Main result]
\label{thm : main}
\begin{align}
\nonumber \frac{1} {2^nn!}&\sum_{\la,\mu \vdash n} {b_{\la,\mu}^n} p_{\la}(x)p_{\mu}(y) =  \sum_{\lambda, \mu \vdash n}m_\lambda(x)m_\mu(y)\sum_{\substack{a,b\\ \{a^1_i,b^1_i,c^1_i\} \in C_{a,b}^{\la}\\ \{a^2_i,b^2_i,c^2_i\} \in C_{a,b}^{\mu}}}r_n(a,b)\times\\
&\prod_{\substack{1\leq i \leq n\\1\leq k\leq 2}}\GenBin{\overline{a}^k_{i-1}-\overline{b}^k_{i-1}}{a^k_i}\VarGenBin{\overline{a}^k_{i-1}-\overline{b}^k_{i}}{b^k_i}{R(\overline{a}^k_{i},\overline{a}^k_{i-1},\overline{b}^k_{i},\overline{b}^k_{i-1},\overline{c}^k_{i-1})}^{c^k_i}
\end{align}
\noindent where the sum runs over pairs of integers $(a,b)$ such that $(a,b,1^{n-a-b})$ is a partition and $C_{a,b}^{\rho}$ is the set of non negative integer sequences $\{a_i,b_i,c_i\}_{1\leq i\leq n}$ with (assume $\rho_i = 0$ for $i>\ell(\rho)$) :
\footnotesize
\begin{flalign*}
&\sum_i a_i = a,\;\;\;\;\sum_i b_i = b,\;\;\;\; a_i+b_i \begin{cases} \in \{\rho_i, \rho_i-1\}&\mbox{ for }i<\ell(\rho)-1\\= \rho_i&\mbox{ otherwise}\end{cases},\;\;\;\;  c_i = \rho_i-a_i-b_i
 \end{flalign*} 
 \normalsize
 Moreover, $\sum_{i=1}^{j} b_i < b$  if  $\sum_{i=1}^{j-1} c_i < c=n-a-b$ and we noted  $\overline{x}^k_i=x-\sum_{j=1}^{i}x^k_i\; (x\in\{a,b,c\})$.
\end{thm}

The formula of theorem \ref{thm : main} doesn't have some of the nice properties of the one in \cite{MV11}. It is not obvious that the coefficients are integers (the summands are rational numbers) and the $r_n(a,b)$ have alternate signs. However, the summation runs over much less parameters and the computation of the coefficients is more efficient especially (see section $\ref{sec : results}$) when $\lambda$ or $\mu$ has a small number of parts.

\section{Connection coefficients of the class algebra}
\label{sec : class}
Let $\chi^\la$ be the irreducible character of $S_n$ indexed by $\la$ and $\chi^\la_\mu$ its value at any element of $C_\mu$. Denote by $f^\la$ the degree $\chi^\la_{1^n}$ of $\chi^\la$. As Biane in \cite{B04} we start with the expression:
\begin{equation}
c^\nu_{\la,\mu}= \frac{n!}{z_\la z_\mu}\sum_{\alpha \vdash n}\frac{\chi^\alpha_\la\chi^\alpha_\mu\chi^\alpha_\nu}{f^\al}
\end{equation}
which becomes much simpler when $\nu=(n)$ as $\chi^\alpha_{(n)} = (-1)^{a}$ if $\alpha = (n-a,1^a)$ and $0$ otherwise. Furthermore $f^{(n-a,1^a)}=\binom{n-1}{a}$. We have:
\begin{equation}
\label{eq : con}
c^n_{\la,\mu} = \frac{n}{z_\la z_\mu}\sum_{a=0}^{n-1}(-1)^a(n-1-a)!a!\chi^{(n-a,1^a)}_\la\chi^{(n-a,1^a)}_\mu
\end{equation}
Then, following \cite{B04}, the generating series are equal to:
\begin{align}
\nonumber \frac{1}{n}&\sum_{\la,\mu \vdash n}c^n_{\la,\mu}p_\la(x)p_\mu(y) =\\
 &\sum_{a=0}^{n-1}(-1)^a(n-1-a)!a!\sum_{\la,\mu  \vdash n}z_\la^{-1}\chi^{(n-a,1^a)}_\la p_\la(x) z_\mu^{-1}\chi^{(n-a,1^a)}_\mu p_\mu(y)
\end{align}
and simplify as $s_\la = \sum_\mu z_\mu^{-1}\chi^{\la}_\mu p_\mu$: (see e.g. \cite{M}):
\begin{equation}
\frac{1}{n}\sum_{\la,\mu \vdash n}c^n_{\la,\mu}p_\la(x)p_\mu(y) = \sum_{a=0}^{n-1}(-1)^a(n-1-a)!a!\sum_{\la,\mu  \vdash n}s_{(n-a,1^a)}(x)s_{(n-a,1^a)}(y)
\end{equation}
In order to recover theorem \ref{thm : mv09}, we need to express the Schur functions in terms of monomials. The transition between the two basis is performed thanks to the Kostka numbers: $s_\la = \sum_\mu K_{\la,\mu} m_\mu$. But obviously (see \cite[I.6. ex 1]{M}) $K_{(n-a,1^a)\la} = \binom{\ell(\la)-1}{a}$:
\begin{align}
\label{eq : bin}
\nonumber &\frac{1}{n}\sum_{\la,\mu \vdash n}c^n_{\la,\mu}p_\la(x)p_\mu(y) =\\
\nonumber  &\sum_{\la,\mu  \vdash n}\sum_{a=0}^{n-1}(-1)^a(n-1-a)!a!\binom{\ell(\la)-1}{a}\binom{\ell(\mu)-1}{a}m_\la(x)m_\mu(y)=\\
& \sum_{\la,\mu  \vdash n}(\ell(\la)-1)!(n-\ell(\la))!\sum_{a=0}^{n-1}(-1)^a\binom{n-1-a}{\ell(\la)-1-a}\binom{\ell(\mu)-1}{a}m_\la(x)m_\mu(y)
\end{align}
Finally, using the classical formulas for binomial coefficients summation (see e.g. \cite{K}), we obtain :
\begin{align}
\nonumber \frac{1}{n}\sum_{\la,\mu \vdash n}c^n_{\la,\mu}p_\la(x)p_\mu(y) &= \sum_{\la,\mu  \vdash n}(\ell(\la)-1)!(n-\ell(\la))!\binom{n-\ell(\mu)}{\ell(\la)-1}m_\la(x)m_\mu(y)\\
&= \sum_{\la,\mu  \vdash n}\frac{(n-\ell(\la))!(n-\ell(\mu))!}{(n+1-\ell(\la)-\ell(\mu))!}m_\la(x)m_\mu(y)
\end{align}
\begin{rem}
The coefficient of $m_\la(x)m_n(y)$ in the series is $n!$. This can be shown directly as $[m_n]p_\mu = 1$ and $\sum_\mu c^n_{\la \mu} = |C_\la| = n!/z_\la$:
\begin{equation}
[m_n(y)]\left(\sum_{\la,\mu \vdash n}c^n_{\la,\mu}p_\la(x)p_\mu(y)\right)= n!\sum_{\la \vdash n}z_\la^{-1}p_\la(x) =n!s_n(x)=n! \sum_{\la \vdash n}m_\la(x)
\end{equation}
\end{rem}

\indent We can recover the formula of theorem \ref{thm : fv10} in a very similar fashion. Equation \ref{eq : con} becomes:
\begin{equation}
c^n_{\la,n} = \frac{1}{z_\la}\sum_{a=0}^{n-1}(n-1-a)!a!\chi^{(n-a,1^a)}_\la
\end{equation}
Pursuing with an analog development gives:
\begin{align}
\sum_{\la \vdash n}c^n_{\la,n}p_\la(x) &= \sum_{\la  \vdash n}(\ell(\la)-1)!(n-\ell(\la))!\sum_{a=0}^{n-1}\binom{n-1-a}{\ell(\la)-1-a}m_\la(x)\\
&= \sum_{\la  \vdash n}(\ell(\la)-1)!(n-\ell(\la))!\binom{n}{n-\ell(\la)+1}m_\la(x)\\
&=\sum_{\la  \vdash n}\frac{n!}{n-\ell(\la)+1}m_\la(x)
\end{align}

\section{Connection coefficients of the double coset algebra}
\label{sec : doublecoset}
Given a partition $\lambda$ and a box $s$ in the Young diagram of $\lambda$, let $l'(s),l(s),$ $a(s),a'(s)$ be the number of boxes to the north, south, east, west of $s$ respectively. These statistics are called {\bf co-leglength, leglength, armlength, co-armlength} respectively. We note as well:
\begin{align}
c_{\lambda}= \prod_{s\in \lambda} (2a(s) + l(s) + 1), \;\;\;\;\;\;\;  c'_{\lambda}= \prod_{s\in \lambda} (2(1+a(s)) + l(s))
\end{align}
Let $\varphi^{\la}_\mu=\sum_{w\in K_{\mu}} \chi^{2\la}_w$ with $2\la=(2\la_1,2\la_2,\ldots,2\la_p)$. Then $| K_\mu |^{-1} \varphi^{\la}_\mu$ is the value of the zonal spherical function indexed by $\lambda$ of the Gelfand pair $(S_{2n}, B_n)$ at the elements of the double coset $K_\mu$. We define as well, $H_{2\la}$ as the product of the hook lengths of the partition $2\la$. We have $H_{2\la} = c_{\lambda}c'_{\lambda}$  \cite[VII.2 eq. (2.19)]{M}. As shown in \cite{HSS92}, the connection coefficients of the double coset algebra verify the relation:
\begin{equation} \label{conncoeff}
b_{\la,\mu}^{\nu} = \frac{1}{|K_{\nu}|} \sum_{\beta \vdash n} \frac{\varphi^{\beta}_\lambda \varphi^{\beta}_\mu\varphi^{\beta}_\nu}{H_{2\beta}} .
\end{equation}
Getting back to the generating series yields:
\begin{equation} 
\sum_{\la,\mu \vdash n} {b_{\la,\mu}^\nu} p_{\la}(x)p_{\mu}(y) = \frac{|B_n|^2}{|K_{\nu}|} \sum_{\beta \vdash n}\frac{\varphi^{\beta}_\nu}{H_{2\beta}}Z_{\beta}(x)Z_{\beta}(y) .
\end{equation}
where $Z_{\beta}(x) = |B_n|^{-1}\sum_\la\varphi^{\beta}_\la p_{\la}(x)$ are the {\bf zonal polynomials} \cite[ch. VII]{M}. As an alternative definition, the $Z_{\beta}(x)$ are special cases of the integral form of the Jack symmetric functions (or polynomials) $J_\beta(x,\alpha)$ with parameter $\alpha = 2$. As described in \cite[VI. ch. 10]{M}, there exist two other classical normalizations of the Jack polynomials. Using notations consistent to \cite{M}, we define $Q_{\lambda}=Z_{\lambda}/c'_{\lambda}$, and $P_{\lambda}=Z_{\lambda}/c_{\lambda}$. The generating series simplifies in the case $\nu = (n)$ since \cite[VII.2 ex 2(c)]{M}:
\begin{equation} \label{macex}
\varphi^{\lambda}_{(n)} = \frac{|K_{(n)}|}{|B_{n-1}|}\prod_{s\in \lambda} (2a'(s)-l'(s))
\end{equation}
where the product omits the square $(1,1)$. As such, $\varphi^{\la}_{(n)}=0$ if $\la \supset (2^3)$ {\em i.e.} if $\la$ is not a near-hook of the form $(a,b,1^{n-a-b})$. Finally, the generating series reduce to:
\begin{align} 
\label{eq gen series}
\nonumber \frac{1}{|B_n|}\sum_{\la,\mu \vdash n} &{b_{\la,\mu}^n} p_{\la}(x)p_{\mu}(y) =\\
 &\frac{|B_n|}{|K_{(n)}|} \sum_{a,b}\varphi^{(a,b,1^{n-a-b})}_{(n)}P_{(a,b,1^{n-a-b})}(x)Q_{(a,b,1^{n-a-b})}(y)
\end{align}
The proof of theorem \ref{thm : main} relies on the evaluation of the expansion in the monomial basis of the zonal polynomials indexed by near hooks, which is done in the following section.

\section{Zonal polynomials on near hooks}
\label{sec : zon poly}
Zonal polynomials are special cases of Jack symmetric functions which in their turn are special cases of the MacDonald polynomials \cite[VI]{M}. While these latter symmetric functions have been heavily studied over the past years (see e.g. \cite{LS03}) no simple expansion in terms monomial functions is known except in the special case of a single part partition (see Stanley \cite{S}). We show that such an expression can be found in the case of Zonal polynomials indexed by near hooks. We start with the general combinatorial formula of MacDonald \cite[VI. eq. (7.13) and (7.13')]{M}:
\begin{eqnarray}
Q_\lambda = \sum_{\mu}m_\mu\sum_{\substack{shape(T)=\lambda,\\ type(T)=\mu}}\phi_{T}, \;\;\;\;\;\;\;\;\;\; P_\lambda = \sum_{\mu}m_\mu\sum_{\substack{shape(T)=\lambda,\\ type(T)=\mu}}\psi_{T}
\end{eqnarray}
\noindent where the internal sums run over all (column strict) tableaux $T$. As usual, a tableau of shape $\lambda$ and type $\mu$ with $\ell(\mu) = p$ can be seen as a sequence of partitions $(\lambda^{(1)},\ldots,\lambda^{(p)})$ such that $\lambda^{(1)}\subset\lambda^{(2)}\subset\ldots\subset\lambda^{(p)}=\lambda$ and each $\lambda^{(i+1)}/\lambda^{(i)}$ (as well as $\lambda^{(1)}$) is a horizontal strip with $|\lambda^{(i+1)}/\lambda^{(i)}| = \mu_{p-i}$ boxes filled with integer $i+1$. 
Further we will assume $\la^{(0)}$ to be the empty partition to simplify notations.

\begin{exm}
\label{ex : tableau}
The following tableau has shape $\lambda = (6,3,1,1)$, type $\mu = (3,2,2,2,2)$ and a sequence $\lambda^{(5)}= \lambda$, $\lambda^{(4)} = (4,3,1)$, $\lambda^{(3)} = (3,2,1)$, $\lambda^{(2)} = (3,1)$, $\lambda^{(1)} = (2)$.  
$$
\vartableau{ 1 & 1 & 2 & 4 & 5 &5\\
2 & 3 & 4\\
3\\
5}
$$
\end{exm}

\begin{figure}[h]
\begin{center}
\includegraphics[width=4cm]{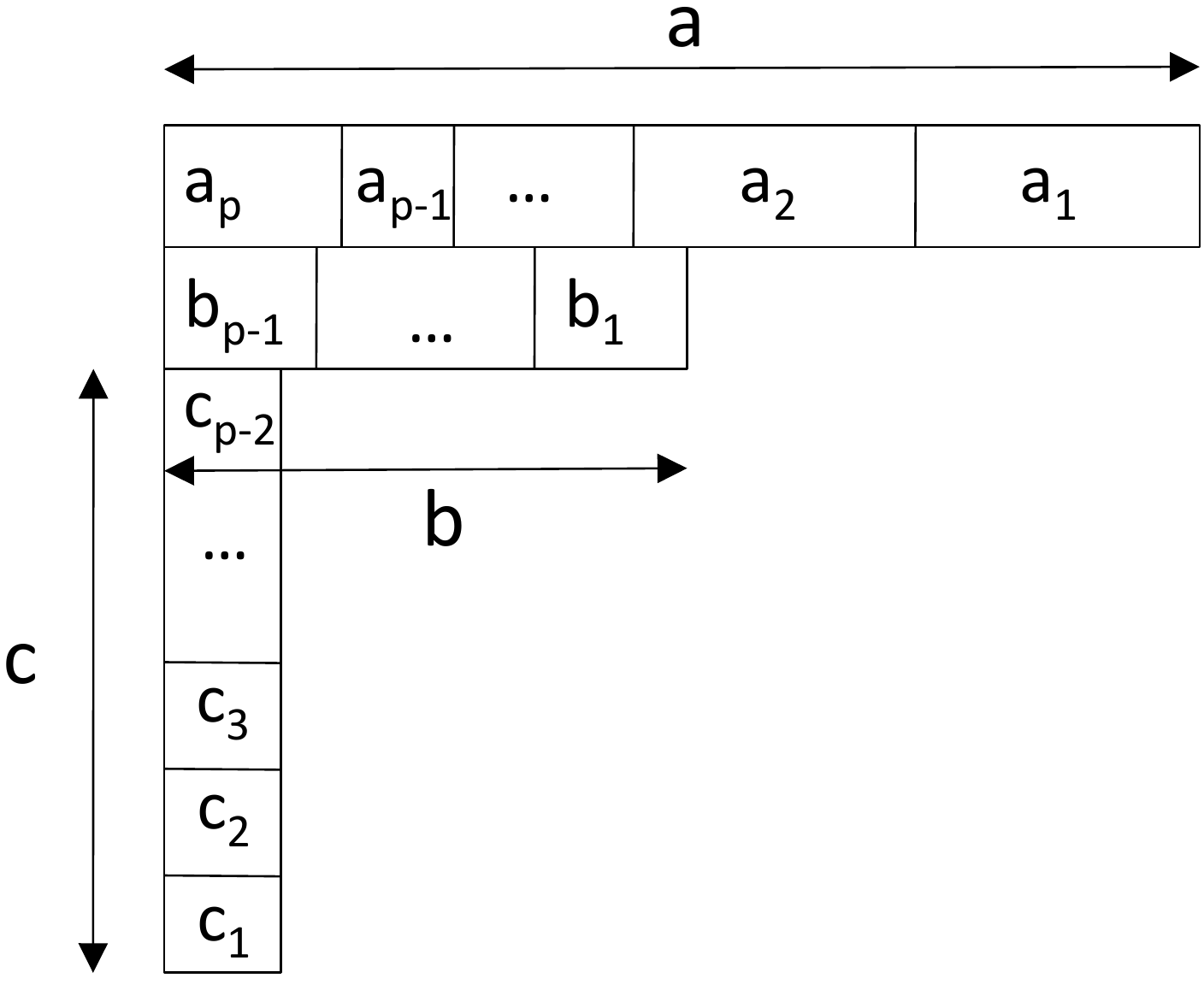}
\end{center}
\caption{Definition of the $(a_i,b_i,c_i)$ on a tableau of shape $(a,b,1^c)$}
\label{tableau}
\end{figure}

\noindent Following (7.11) of \cite{M}:
\begin{equation}
\phi_{T} = \prod_{i=0}^{p-1}\phi_{\lambda^{(i+1)}/\lambda^{(i)}}, \;\;\;\;\;\;\; \psi_{T} = \prod_{i=0}^{p-1}\psi_{\lambda^{(i+1)}/\lambda^{(i)}}
\end{equation}
The analytic formulation of $\phi_{\lambda/\mu}$ and $\psi_{\lambda/\mu}$ is given by \cite[VI. 7. ex 2.(a)]{M} for the general case of MacDonald polynomials. We use the formulation given by Okounkov and Olshanski in \cite{OO97} (equation (6.2)) specific to the Jack symmetric functions:
\begin{equation}
\psi_{\lambda/\mu} = \prod_{1\leq k \leq l \leq \ell(\mu)}\frac{(\mu_k-\la_{k+1}+\theta(l-k)+1)_{\theta-1}(\la_k-\mu_{k}+\theta(l-k)+1)_{\theta-1}}{(\la_k-\la_{l+1}+\theta(l-k)+1)_{\theta-1}(\mu_k-\mu_{l}+\theta(l-k)+1)_{\theta-1}}
\end{equation}
where $(t)_r=\Gamma(t+r)/\Gamma(t)$ for any numbers $t$ and $r$. This formula holds in the case of zonal polynomials by setting $\theta=1/2$.
Similarly,
\begin{equation}
\label{formula : phi}
\phi_{\lambda/\mu} = \prod_{1\leq k \leq l \leq \ell(\la)}\frac{(\mu_k-\la_{k+1}+\theta(l-k)+1)_{\theta-1}(\la_k-\mu_{k}+\theta(l-k)+1)_{\theta-1}}{(\la_k-\la_{l}+\theta(l-k)+1)_{\theta-1}(\mu_k-\mu_{l+1}+\theta(l-k)+1)_{\theta-1}}
\end{equation}
Given a tableau of shape $(a,b,1^c)$ of type $\mu = (\mu_1,\mu_{2},\ldots, \mu_p)$ we note $a_i$ (resp. $b_i$ and $c_i$) the number of boxes of the first line (resp. second line and column) filled with integer $p-i+1$. This notation is illustrated by figure \ref{tableau}.

\noindent Obviously, using the notations of section \ref{sec : intro} we have :
\begin{align*}
&\lambda^{(p-i)} = (\overline{a}_i,\overline{b}_i,1^{\overline{c}_i}),  \;\;\;\;\;\;  a_i+b_i+c_i = \mu_i  \\
&a_i\leq \overline{a}_{i-1}-\overline{b}_{i-1}, \;\;\;\;\;\;\;\;\;\;\; b_i\leq \overline{b}_{i-1}-1 \mbox{      if       } \overline{c}_{i} > 0\\
&b_p = c_p = c_{p-1} = 0, \;\;\;\;\;\;\;  c_i \in \{0,1\}
\end{align*}

\begin{exm}
The tableau of example \ref{ex : tableau} is of shape $(6,3,1^2)$ and check $(a_1, b_1$, $c_1)$ $= (2,0,1)$, $(a_2, b_2, c_2) = (1,1,0)$, $(a_3, b_3, c_3) = (0,1,1)$, $(a_4, b_4, c_4) = (1,1,0)$,  $(a_5, b_5, c_5) = (2,0,0)$.
\end{exm}
Given a tableau of shape $(a,b,1^c)$ with filling described by $(a_i,b_i,c_i)_{1\leq i\leq p}$ we apply formula \ref{formula : phi} to $\phi_{\lambda^{(p-i+1)}/\lambda^{(p-i)}}$. The following are the non-$1$ contributions to the product:

\begin{itemize}
\item[(i)] for $k=l=1$, the factor in equation \ref{formula : phi} is $\frac{(a_i+1)_{-1/2}(\overline{a}_{i-1}-\overline{b}_{i-1}-a_i+1)_{-1/2}}{(1)_{-1/2}(\overline{a}_{i-1}-\overline{b}_{i-1}-a_i+b_i+1)_{-1/2}}$
\item[(ii)] for $k=1$, $l=2$: $\frac{(\overline{a}_{i-1}-\overline{b}_{i-1}+b_i+3/2)_{-1/2}}{(\overline{a}_{i-1}-\overline{b}_{i-1}+3/2)_{-1/2}}$ 
\item[(iii)] for $k=l=2$: $\frac{(b_i+1)_{-1/2}}{(1)_{-1/2}}$
\end{itemize}
\noindent Additionally, if $c_i=1$ we have:
\begin{itemize}
\item[(iv)]  $k=1$, $l =\overline{c}_{i-1}+1$: $\frac{(\overline{a}_{i-1}-a_i+\overline{c}_{i-1}/2)_{-1/2}}{(\overline{a}_{i-1}-a_i+\overline{c}_{i-1}/2+1)_{-1/2}}$   
\item[(v)] $k=1$, $l =\overline{c}_{i-1}+2$: $\frac{(\overline{a}_{i-1}+(\overline{c}_{i-1}+1)/2+1)_{-1/2}}{(\overline{a}_{i-1}+(\overline{c}_{i-1}+1)/2)_{-1/2}}$  
\item[(vi)] $k=2$, $l =\overline{c}_{i-1}+1$: $\frac{(\overline{b}_{i-1}-b_i+(\overline{c}_{i-1}-1)/2)_{-1/2}}{(\overline{b}_{i-1}-b_i+(\overline{c}_{i-1}-1)/2+1)_{-1/2}}$
\item[(vii)] $k=2$, $l =\overline{c}_{i-1}+2$: $\frac{(\overline{b}_{i-1}+\overline{c}_{i-1})/2+1)_{-1/2}}{(\overline{b}_{i-1}+\overline{c}_{i-1}/2)_{-1/2}}$  
\item[(viii)] $k>2$, $l =\overline{c}_{i-1}+1$: $\frac{(\overline{c}_{i-1}+1-k)/2+1)_{-1/2}}{(\overline{c}_{i-1}+1-k)/2+2)_{-1/2}}$
\item[(ix)] $k>2$, $l =\overline{c}_{i-1}+2$: $\frac{(\overline{c}_{i-1}+2-k)/2+2)_{-1/2}}{(\overline{c}_{i-1}+2-k)/2+1)_{-1/2}}$   
\end{itemize} 
As $(t+1)_{-1/2} = (2t)!/(\sqrt{\pi}t!^22^{2t})$ and $(t+3/2)_{-1/2} = \sqrt{\pi}t!^22^{2t+1}/(2t+1)!$, (i)$\times$(ii)$\times$(iii) gives:
\begin{equation}
\frac{1}{2^{2a_i-2b_i}}\GenBin{a_{i-1}-b_{i-1}}{a_i}\VarGenBin{a_{i-1}-b_{i}}{b_i}\frac{\binom{2a_{i-1}-2b_{i-1}}{a_{i-1}-b_{i-1}}}{\binom{2a_{i}-2b_{i}}{a_{i}-b_{i}}}
\end{equation}
\noindent Then, as $(t)_{-1/2}/(t+1)_{-1/2} = 2t/(2t-1)$, combining (iv), (v), (vi) and (vii) yields $1$ if $c_i=0$ and
\begin{equation}
\label{eq : R}
\frac{2\overline{a}_{i}+\overline{c}_{i-1}}{2\overline{a}_{i}+\overline{c}_{i-1}-1}\frac{2\overline{a}_{i-1}+\overline{c}_{i-1}}{2\overline{a}_{i-1}+\overline{c}_{i-1}+1}\frac{2\overline{b}_{i}+\overline{c}_{i-1}-1}{2\overline{b}_{i}+\overline{c}_{i-1}-2}\frac{2\overline{b}_{i-1}+\overline{c}_{i-1}-1}{2\overline{b}_{i-1}+\overline{c}_{i-1}}
\end{equation}
\noindent otherwise. For all $k>2$ the combination of (viii) and (ix) can be written as
\begin{equation}
\prod_{k>2}^{\overline{c}_{i-1}+1}\frac{\overline{c}_{i-1}+3-k}{\overline{c}_{i-1}+3-(k+1)}\prod_{k>2}^{\overline{c}_{i-1}+2}\frac{\overline{c}_{i-1}+4-(k+1)}{\overline{c}_{i-1}+4-k}=\frac{\overline{c}_{i-1}}{\overline{c}_{i-1}+1} =\frac{\overline{c}_{i}+1}{\overline{c}_{i-1}+1}
\end{equation}
\noindent The last ratio holds in the general case ($c_i \in \{0,1\}$). Putting everything together and multiplying the factors for $i=1\ldots p$ gives the formula:

\begin{thm}[Zonal polynomials on near hooks]
\label{thm : Zonal}
\begin{align}
\nonumber &Q_{a,b,1^c}= \frac{\binom{2a-2b}{a-b}}{4^{a-b}(1+c)}\times\\
&\sum_{\substack{\lambda \vdash a+b+c\\\{a_i,b_i,c_i\}}}\prod_{i}\GenBin{\overline{a}_{i-1}-\overline{b}_{i-1}}{a_i}\VarGenBin{\overline{a}_{i-1}-\overline{b}_{i}}{b_i}{R(\overline{a}_{i},\overline{a}_{i-1},\overline{b}_{i},\overline{b}_{i-1},\overline{c}_{i-1})}^{c_i}
\end{align}
\noindent Similarly, we find:
\begin{align}
\nonumber &P_{a,b,1^c}=  f(a,b,c)\times\\
& \sum_{\substack{\lambda \vdash a+b+c\\\{a_i,b_i,c_i\}}}\prod_{i}\GenBin{\overline{a}_{i-1}-\overline{b}_{i-1}}{a_i}\VarGenBin{\overline{a}_{i-1}-\overline{b}_{i}}{b_i}{R(\overline{a}_{i},\overline{a}_{i-1},\overline{b}_{i},\overline{b}_{i-1},\overline{c}_{i-1})}^{c_i}
\end{align}
\footnotesize
where $f(a,b,c)=\begin{cases} (2a+c+1)(2b+c)/\left ((2a+c)(2b+c-1)\VarGenBin{a-1}{b-1}\right) \mbox{ if } b\neq 0\\  1 \;\;\;\mbox{ otherwise}\end{cases}$
\normalsize
\end{thm}
\begin{exm}
When $b=c=0$, we get  $P_{n}= \sum_{\lambda \vdash n} \GenBin{n}{\lambda}m_\lambda$ and $Q_{n}=\frac{1}{4^{n}}\binom{2n}{n} \sum_{\lambda \vdash n} \GenBin{n}{\lambda}m_\lambda$,
which is equivalent to the formula of Stanley in \cite{S} (the extension of {\footnotesize$\GenBin{x}{y}$} to multinomial coefficients is straightforward) .
\end{exm}
\begin{exm}
Only the following tableau of shape $(a,b,1^{c})$ yields a non zero contribution to the coefficient in $m_{a,b,1^{c}}$ of $P_{a,b,1^{c}}$:
\footnotesize
$$
\tableau{ 1 &p\hspace{-1mm}-\hspace{-1.4mm}1& \ldots &p\hspace{-1mm}-\hspace{-1.4mm}1&p\hspace{-1mm}-\hspace{-1.4mm}1&p& \ldots&p & p\\
2 &p&\ldots & p & p\\
3\\
:\\
p\hspace{-1mm}-\hspace{-1.4mm}2\\
p\hspace{-1mm}-\hspace{-1.4mm}1\\
p}
$$
\normalsize
where we have $(a-b)$ cases filled with $p$ on the first line, $(b-1)$ on the second and $1$ at the bottom of the column. Then $(b-1)$ cases filled with $p-1$ on the first line and $1$ in the column. Finally the column is labeled from bottom to top with $p-2,p-3,\ldots,1$. When all the $c_i$'s are equal to $1$ expression \ref{eq : R} reads:
\vspace{-1mm}
\begin{equation}
\frac{2\overline{a}_{i}+\overline{c}_{i}+1}{2\overline{a}_{i}+\overline{c}_{i}}\frac{2\overline{a}_{i-1}+\overline{c}_{i-1}}{2\overline{a}_{i-1}+\overline{c}_{i-1}+1}\frac{2\overline{b}_{i}+\overline{c}_{i}}{2\overline{b}_{i}+\overline{c}_{i}-1}\frac{2\overline{b}_{i-1}+\overline{c}_{i-1}-1}{2\overline{b}_{i-1}+\overline{c}_{i-1}}
\end{equation}
\noindent and we have 
\begin{equation}
\prod_i{R(\overline{a}_{i},\overline{a}_{i-1},\overline{b}_{i},\overline{b}_{i-1},\overline{c}_{i-1})}^{c_i} = \frac{3}{2}\times\frac{2a+c}{2a+c+1}\times\frac{2}{1}\times\frac{2b+c-1}{2b+c}
\end{equation}
The other contributing factors read
\begin{align}
\nonumber \GenBin{a-b}{a-b}\VarGenBin{a-1}{b-1}\GenBin{b-1}{b-1}&\VarGenBin{b-1}{0}\GenBin{0}{0}\VarGenBin{1}{1}\GenBin{1}{1}\VarGenBin{1}{0} \\
 &=\frac{1}{3}\VarGenBin{a-1}{b-1}
\end{align}
where we use the fact that $\VarGenBin{1}{1} = 1/3$. Putting everything together, yields the classical property:
\begin{equation}
[m_{a,b,1^{c}}]P_{a,b,1^{c}}=1
\end{equation}
\end{exm}

\section{Proof of main theorem, examples and further results}
\label{sec : results}
The proof  of theorem \ref{thm : main} follows immediately from equation \ref{eq gen series}, theorem \ref{thm : Zonal} and the final remark:
\begin{equation}
\prod_{s\in (a,b,1^c)} (2a'(s)-l'(s)) = \begin{cases}(-1)^{c+1}(c+1)!(2a-2)!!(2b-3)!!\mbox{ if } b>0\\ (2a-2)!! \mbox{ otherwise}\end{cases}
\end{equation}
Now we study some particular coefficients of the generating series in theorem \ref{thm : main} and emphasize the link with the alternative formula in \cite{MV11}.\\
\subsection{Coefficient of $m_\la(x)m_n(y)$}
As a first example, we notice that only the near hooks $(a,b,1^c)$ that have less parts than $\min(\ell(\lambda),\ell(\mu))$ contribute to the coefficient of $m_\la(x)m_\mu(y)$ in theorem \ref{thm : main}. If either $\la$ or $\mu$ is the single part partition $(n)$, then only the one row tableaux of length $n$ contribute to the coefficient in the generating series. It is straightforward from theorem \ref{thm : main} that
\begin{equation}
[m_n(x)m_n(y)]\left (\frac{1} {2^nn!}\sum_{\la,\mu \vdash n} {b_{\la,\mu}^n} p_{\la}(x)p_{\mu}(y)\right ) = r_n(n,0) = (2n-1)!!
\end{equation}
From the perspective of the combinatorial interpretation in \cite{MV11} this is an obvious result as the coefficient of $m_n(x)m_n(y)$ is the number of pairings on a set of size $2n$. More interestingly, noticing that
\begin{equation}
\GenBin{n}{\lambda} = \binom{n}{\la}\frac{(2\la-1)!!}{(2n-1)!!}
\end{equation}
\noindent where $(2\la-1)!! = \prod_i(2\la_i-1)!!$, we get:
\begin{equation}
\label{eq : mlamn}
[m_\la(x)m_n(y)]\left (\frac{1} {2^nn!}\sum_{\la,\mu \vdash n} {b_{\la,\mu}^n} p_{\la}(x)p_{\mu}(y)\right ) =  \binom{n}{\la}(2\la-1)!!
\end{equation}
This result is not obvious to derive from the formula in \cite{MV11} obtained by Lagrange inversion, but it can be proved with the combinatorial interpretation in terms of some decorated bicolored forests with a single white vertex. 
The exact definition of these forests is given in \cite[Def. 2.10]{MV11}. We briefly remind that they are composed of white and black (internal or root) vertices.
The descendants of a given vertex are composed of edges (linking a white vertex and a black one), thorns and loops. Additionally, there is a bijection between thorns connected to white vertices and thorns connected to black vertices and a mapping of the loops on white (resp. black) vertices to the set of black (resp. white) vertices. In this paper, the authors show that $(2n-1)!! =F_n$ is the number of two-vertex (one white and one black) forests of size $n$ (see figure \ref{forests} for examples). 
\begin{figure}[h]
\begin{center}
\includegraphics[width=6cm]{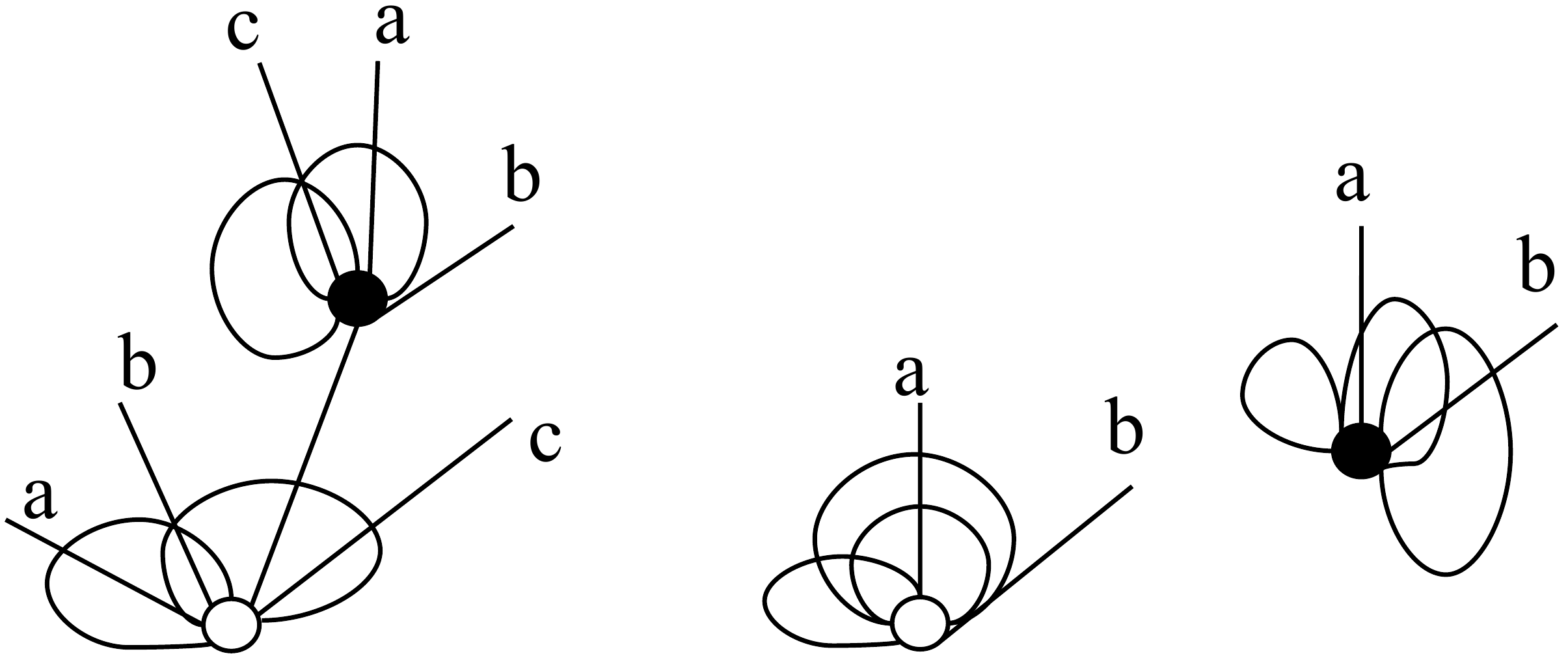}
\end{center}
\caption{Two examples of two-vertex forests for $n=7$ (left) and $8$ (right). Letters depict the bijection between thorns}
\label{forests}
\end{figure}
But a forest with one white (root) vertex and $\ell(\lambda)$ black vertices of degree distribution $\la$ ($F_\la$ denotes the number of such forests) can be seen as a $\ell(\lambda)$-tuple of two vertex forests of size $\la_i$. The $i$-th forest is composed of the $i$-th black vertex with its descendants and one white vertex with a subset of descendants of the original one's containing (i) the edge linking the white vertex and the $i$-th black vertex, (ii) the thorns in bijection with the thorns of the $i$-th black vertex, (iii) the loops mapped to the $i$-th black vertex. 
The construction is bijective if we distinguish in the initial forest the black vertices with the same degree ($Aut_\la$ ways to do it) and we keep track in the tuple of forests the initial positions of the descendants of the white vertices within the initial forest ($\binom{n}{\la}$ possible choices). We get:
\begin{equation}
Aut_\la F_\la = \binom{n}{\la}\prod_i F_{\la_i} =\binom{n}{\la}(2\la-1)!!
\end{equation}
\begin{figure}[h]
\begin{center}
\includegraphics[width=10cm]{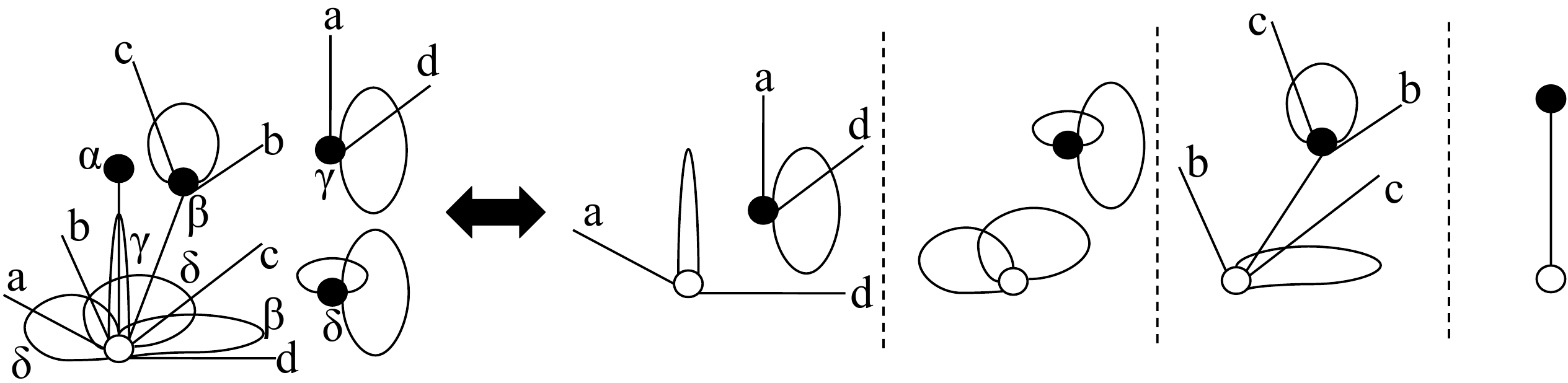}
\end{center}
\caption{Splitting a forest of black degree distribution $\la$ into a $\ell(\la)$-tuple of two vertex forests for $\la = (4^3,1)$. Greek letters depict the mapping on the sets of loops connected to the white vertex}
\label{forests}
\end{figure}
Finally, according to \cite{MV11}, $Aut_\la F_\la$ is equal to the desired coefficient in (\ref{eq : mlamn}). 
\begin{rem}
Using the main formula of \cite{MV11} we've shown:
\begin{equation}
\sum_{Q,Q'}\prod_{i,j}\frac{2^{Q'_{ij}-2j(Q_{ij}+Q'_{ij})}}{Q_{ij}!Q'_{ij}!}{\binom{i-1}{j,j}}^{Q_{ij}}{\binom{i-1}{j,j-1}}^{Q'_{ij}} = \frac{(2\la-1)!!}{\la!Aut_\la}
\end{equation}
where the sum runs over matrices $Q$ and $Q'$ with $m_i(\lambda)=\sum_{j \geq 0}Q_{ij} + Q'_{ij}$.
\end{rem}
\begin{rem}
If we admit the expansion of $Z_n$ in terms of monomials, we can directly show (\ref{eq : mlamn}) as: 
\begin{align}
  \nonumber {|B_n|}^{-1}\sum_{\la \vdash n}\left ( \sum_{\mu \vdash n}{b_{\la,\mu}^n}\right) p_{\la}(x)&=  {|B_n|}^{-1}\sum_{\la \vdash n}|K_\la|p_{\la}(x) =Z_n(x)\\
  &= \sum_{\la \vdash n}\binom{n}{\la}(2\la-1)!!m_\la(x)
\end{align}

\end{rem}
\subsection{Coefficient of $m_{n-p,1^{p}}(x)m_{n-p,1^{p}}(y)$}
The number $F_{(n-p,1^p),(n-p,1^p)}$ of forests with $p+1$ white and $p+1$ black vertices, both of degree distribution $(n-p,1^{p})$, can be easily obtained from the number of two-vertex forests $F_{n-2p}$. We consider $2p\leq n-1$, it's easy to show the coefficient to be equal to $0$ otherwise. Two cases occur: either the white vertex with degree $n-p$ is the root and there are $\binom{n-p}{p}\times\binom{n-p-1}{p}$ ways to add the black and the white descendants of degree $1$, or the root is a white vertex of degree $1$ and there are $\binom{n-p-1}{p-1}\times\binom{n-p-1}{p}$ ways to add the remaining white vertices and the $p$ black vertices of degree $1$. We have:
\begin{equation}
F_{(n-p,1^p),(n-p,1^p)} = F_{n-2p}\binom{n-p-1}{p}\left [\binom{n-p}{p}+\binom{n-p-1}{p-1}\right ]
\end{equation}
As a result, we obtain :
\begin{align}
\nonumber [m_{n-p,1^{p}}(x)m_{n-p,1^{p}}(y)]&\left (\frac{1} {2^nn!}\sum_{\la,\mu \vdash n} {b_{\la,\mu}^n} p_{\la}(x)p_{\mu}(y)\right )\\
&\nonumber = Aut_{n-p,1^{p}}^2F_{(n-p,1^p),(n-p,1^p)}\\
&= n(n-2p)\left(\frac{(n-p-1)!}{(n-2p)!}\right)^2(2n-4p-1)!!
\end{align}
We check this result with the formula of theorem \ref{thm : main} for the special cases $p\in\{1,2\}$. For $p=1$, two tableaux 
\begin{equation}
\overbrace{\vartableau{1 & 2&... & 2&2 }}^{n}\;\;\mbox{ and }\;\;\overbrace{\vartableau{1 & 2&... & 2&2\\2 }}^{n-1}
\end{equation}
contribute to the coefficient with respective contributions $n^2(2n-5)!!(2n-3)/(2n-1)$ and $-2n(2n-5)!!(n-1)/(2n-1)$. Adding them gives the desired result $n(n-2)(2n-5)!!$. In the case $p=2$, the following tableaux are contributing: 
\begin{equation}
\overbrace{\vartableau{1&2&3&...&3&3}}^{n},\;\;\;\;  \overbrace{\vartableau{1&2&3&...&3\\3}}^{n-1},\;\;\;\;  \overbrace{\vartableau{1&3&3&...&3\\2}}^{n-1},\;\;\;\;  \overbrace{\vartableau{1&2&3&...&3\\3&3}}^{n-2},\;\;\;\overbrace{\vartableau{1&3&3&...&3\\2\\3}}^{n-2}
\end{equation}
The second and third tableaux are the two possible fillings for the shape $a=n-1$ and $b=1$ and we have to consider the cross contributions (one filling for $\la$ and the other for $\mu$). Combining thecontributions for tableaus 2 and 3, the sum writes :
\footnotesize
\begin{align}
\nonumber &\frac{n^2(n-1)^2(2n-9)!!(2n-7)(2n-5)}{(2n-1)(2n-3)}-\frac{2n(2n^2-6n+3)^2(2n-9)!!(2n-7)}{(n-1)(2n-1)(2n-5)}\\
\nonumber&-\frac{8n(n-2)^2(n-3)^2(2n-9)!!}{3(2n-5)(2n-3)}+\frac{2n(2n-7)(2n-3)(2n-9)!!}{3(n-1)}
\end{align}
\normalsize
Finally, we have the desired result: $n(n-4)(n-3)^2(2n-9)!!$.
\normalsize
\subsection{Generating series for $b_{\la n}^n$}
As a final application we compute the generating series for $b_{\la n}^n$ in a similar fashion as in section \ref{sec : doublecoset} and \ref{sec : zon poly}:
\begin{equation} 
\Pi_n = \frac{1}{|B_n|}\sum_{\la \vdash n} {b_{\la n}^n} p_{\la} = \frac{1}{|K_{(n)}|} \sum_{a,b}\frac{\left (\varphi^{(a,b,1^{n-a-b})}_{(n)}\right )^2}{c'_{a,b,1^{n-a-b}}}P_{(a,b,1^{n-a-b})}
\end{equation}
Using the same notations as in theorem \ref{thm : main}, we find:
\begin{thm}[Generating series for $b_{\la n}^n$]
\begin{align}
\nonumber &\Pi_n =  \sum_{\substack{\la \vdash n \\ a,b \\ a_i,b_i,c_i}}r'_n(a,b)\\
&\times \prod_{1\leq i \leq n}\GenBin{\overline{a}_{i-1}-\overline{b}_{i-1}}{a_i}\VarGenBin{\overline{a}_{i-1}-\overline{b}_{i}}{b_i}{R(\overline{a}_{i},\overline{a}_{i-1},\overline{b}_{i},\overline{b}_{i-1},\overline{c}_{i-1})}^{c_i}m_\lambda 
\end{align}
with $r'_n(n,0) = (2n-2)!!$ and $r'_n(x,y) = 2n\frac{(n+1-x-y)!(2x-2)!!(2y-3)!!}{(n+x-y)(n+y-x-1)}\; (y>0)$
\end{thm}
Contrary to $r_n$, $r'_n$ is not of alternate sign. This leaves possibilities for asymptotic evaluations.
\begin{exm}
The following table gives the value of some coefficients in the monomial expansion of $\Pi_n$.\\
\footnotesize
\begin{tabular}{|c|c|c|}
\hline
$\la$&$(n)$&$(n-1,1)$\\
\hline
$[m_\la]\Pi_n$&$(2n-2)!!$&$n(2n-4)!!$\\
\hline
$\la$&$(n-2,1,1)$&$(n-3,1,1,1)$\\
\hline
$[m_\la]\Pi_n$&$n(n-1)(2n-6)!!$&$n(n-1)(n-2)(2n-8)!!$\\
\hline
$\la$&$(1^n)$&$(n-2,2)$\\
\hline
$[m_\la]\Pi_n$&$n!$&$n(2n-6)!!(3n-5)/2$\\ 
\hline
$\la$&$(n-3,3)$&$(n-4,4)$\\
\hline
$[m_\la]\Pi_n$&$n(2n-8)!!(5n^2-21n+20)/2$&$n(2n-10)!!(35n^3-270n^2+649n-486)/8$\\
\hline
\end{tabular}
\end{exm}
\normalsize
\vspace{-2mm}
\section*{Acknowledgements}
\label{sec:ack}
The author acknowledges the support of ERC under "ERC StG 208471 - ExploreMaps" and thanks Alejandro Morales for interesting discussions about the state of the art in the field.


\bibliographystyle{alpha}
\bibliography{biblio}

%

%
%
%

\end{document}